\documentclass{amc}
\usepackage{amssymb,amsmath,float,nicefrac} 
\usepackage{euscript}
\usepackage[dvips]{epsfig}
\usepackage{graphics} 

\def\currentvolume{1}
 \def\currentissue{1}
  \def\currentyear{2007}
   
    \def\ppages{131-149}

\newtheorem{theorem}{Theorem}
\newtheorem{corollary}{Corollary}

\newtheorem{lemma}{Lemma}

\theoremstyle{definition}

\def\thesubsection{\thesection-\Alph{subsection}} 

\def\dist{\qopname\relax{no}{dist}}

\def\cosec{\qopname\relax{no}{cosec}}

\def\Cap{\qopname\relax{no}{Cap}}

\newcommand\nc\newcommand
\nc\bfa{{\mathbf a}}\nc\bfA{{\mathbf A}}\nc\cA{{\mathcal A}}
\nc\bfb{{\mathbf b}}\nc\bfB{{\mathbf B}}\nc\cB{{\mathcal B}}
\nc\bfc{{\mathbf c}}\nc\bfC{{\mathbf C}}\nc\cC{{\mathcal C}}
\nc\bfd{{\mathbf d}}\nc\bfD{{\mathbf D}}\nc\cD{{\mathcal D}}
\nc\bfe{{\mathbf e}}\nc\bfE{{\mathbf E}}\nc\cE{{\mathcal E}}
\nc\bff{{\mathbf f}}\nc\bfF{{\mathbf F}}\nc\cF{{\mathcal F}}
\nc\bfg{{\mathbf g}}\nc\bfG{{\mathbf G}}\nc\cG{{\mathcal G}}
\nc\bfh{{\mathbf h}}\nc\bfH{{\mathbf H}}\nc\cH{{\mathcal H}}
\nc\bfi{{\mathbf i}}\nc\bfI{{\mathbf I}}\nc\cI{{\mathcal I}}
\nc\bfj{{\mathbf j}}\nc\bfJ{{\mathbf J}}\nc\cJ{{\mathcal J}}
\nc\bfk{{\mathbf k}}\nc\bfK{{\mathbf K}}\nc\cK{{\mathcal K}}
\nc\bfl{{\mathbf l}}\nc\bfL{{\mathbf L}}\nc\cL{{\mathcal L}}
\nc\bfm{{\mathbf m}}\nc\bfM{{\mathbf M}}\nc\cM{{\mathcal M}}
\nc\bfn{{\mathbf n}}\nc\bfN{{\mathbf N}}\nc\cN{{\mathcal N}}
\nc\bfo{{\mathbf o}}\nc\bfO{{\mathbf O}}\nc\cO{{\mathcal O}}
\nc\bfp{{\mathbf p}}\nc\bfP{{\mathbf P}}\nc\cP{{\mathcal P}}
\nc\bfq{{\mathbf q}}\nc\bfQ{{\mathbf Q}}\nc\cQ{{\mathcal Q}}
\nc\bfr{{\mathbf r}}\nc\bfR{{\mathbf R}}\nc\cR{{\mathcal R}}
\nc\bfs{{\mathbf s}}\nc\bfS{{\mathbf S}}\nc\cS{{\mathcal S}}
\nc\bft{{\mathbf t}}\nc\bfT{{\mathbf T}}\nc\cT{{\mathcal T}}
\nc\bfu{{\mathbf u}}\nc\bfU{{\mathbf U}}\nc\cU{{\mathcal U}}
\nc\bfv{{\mathbf v}}\nc\bfV{{\mathbf V}}\nc\cV{{\mathcal V}}
\nc\bfw{{\mathbf w}}\nc\bfW{{\mathbf W}}\nc\cW{{\mathcal W}}
\nc\bfx{{\mathbf x}}\nc\bfX{{\mathbf X}}\nc\cX{{\mathcal X}}
\nc\bfy{{\mathbf y}}\nc\bfY{{\mathbf Y}}\nc\cY{{\mathcal Y}}
\nc\bfz{{\mathbf z}}\nc\bfZ{{\mathbf Z}}\nc\cZ{{\mathcal Z}}
\nc\od{{\bar d}}\nc\ow{{\bar w}}\nc\odelta{{\bar\delta}}
\nc\ox{{\bar x}}\nc\oy{{\bar y}}\nc\ou{{\bar u}}
\nc\oh{{\bar h}}
\newcommand{\half}{\nicefrac12}

\newcommand{\remove}[1]{}\newcommand\reals{{\mathbb R}}

\nc\dgv{\delta_{\text{\rm GV}}}
\nc\rcrit{R_{\text{\rm crit}}}
\nc\Esp{E_{\text{\rm sp}}}
\renewcommand\epsilon{\varepsilon}
\renewcommand\angle{\measuredangle}

\title[Codes in spherical caps]	
{Codes in spherical caps} 

\author[A. Barg and O.~R. Musin]{}


\email{abarg@umd.edu,omusin@gmail.com}

\thanks{The first author was supported in part by NSF grants 
CCR0310961, CCF0515124, and by NSA grant H98230-06-1-0044. The second
author was supported in part by NSF grant CCR0310961.}

\begin{document}
\maketitle

\centerline{\scshape Alexander Barg}
\medskip
{\footnotesize
 \centerline{Dept. of ECE and Institute for Systems Research}
  \centerline{University of Maryland}
   \centerline{College Park, MD 20742, USA}
   }

\medskip
\centerline{\scshape Oleg R. Musin}
\medskip
{\footnotesize
 \centerline{Institute for Math. Study of Complex Systems}
  \centerline{Moscow State University}
   \centerline{Moscow, Russia}
   }

\medskip

 \medskip

\begin{abstract}
We consider bounds on codes in spherical caps and related problems
in geometry and coding theory. An extension of the Delsarte method
is presented that relates upper bounds on the size of spherical codes
to upper bounds on codes in caps. Several new upper bounds on 
codes in caps are derived. Applications of these bounds 
to estimates of the kissing numbers and one-sided kissing numbers are
considered.

It is proved that the maximum size of codes in spherical caps for large
dimensions is determined by the maximum size of spherical codes, so 
these problems are asymptotically equivalent.
\end{abstract}

\section{Introduction}

The subject of this paper is codes in spherical caps, i.e., packings
of a metric ball on the surface of the sphere in $\reals^n$ (a spherical cap)
with metric balls (caps) of a smaller radius. 
Codes in spherical caps are related to
more familiar spherical codes and find a number of interesting applications
in both classical and recent works. 

Spherical cap codes
have been used to derive an asymptotic upper bound on the maximum
size of spherical codes and a bound on the packing density $\Delta_n$
of the $n$-space by equal spheres, see Sidelnikov \cite{sid74},
Kabatiansky and Levenshtein \cite{kab78} and Levenshtein \cite{lev83a}.
More recently they have been used to derive upper bounds on the 
size of binary constant weight codes, see 
Agrell, Vardy and Zeger \cite{agr00}.
Even more recently, estimates of the maximum size of codes in a spherical 
cap have been
used together with an extension of Delsarte's method to derive new
estimates of the kissing number $k(n)$ in small dimensions.
In particular, a long-standing
conjecture that $k(4)=24$ was solved in \cite{mus03,mus05b} and
a related problem
for ``one-sided kissing numbers'' was solved in dimension 4 in \cite{mus05c}.

In this paper we focus on the study of spherical cap codes rather
than on their applications in related geometric problems of coding theory.
More specifically, we study bounds on the size of spherical cap
codes with a given angular separation and their relation to spherical
codes. 

In Section \ref{sect:ls} we recall a few known bounds on the size of
codes in a spherical cap and spherical strip. 
In Section \ref{sect:extension} we formulate an extension of Delsarte's
bound on the size of spherical codes to cover the case of spherical caps. 
As usual, the polynomial involved in
the computation of the bound must be expandable in a linear combination
of Gegenbauer polynomials with nonnegative coefficients. A new condition
in the theorem relates the values of the polynomial to constructions
of codes in a spherical cap. The method described
was used in \cite{mus05b} to prove that $k(4)=24.$
This link serves as an additional motivation
for studying spherical cap codes. 

In Section \ref{sect:large} we show that in the case of large code
distance, the size of a code in a spherical cap can be exactly expressed
via the size of codes on the entire sphere. The result is used to relate
the size of spherical cap codes with the kissing number $k(n).$
We consider examples of small dimensions $n=3,4$ and 
illustrate the application of the extended Delsarte's method to the
derivation of the values of $k(n)$ in these cases.

In Section \ref{sect:bound} we derive a new bound on the size of
spherical cap codes that relies on a transformation from codes in caps
to codes on the hemisphere. 
In the same section we also address the problem of the maximum size of
spherical cap codes in the case of large dimensions. 
A common perception in coding theory, originating with the asymptotic
results of \cite{kab78} is that codes in spherical caps are
analogous to constant weight binary codes (i.e., codes formed of vectors with a
fixed number of ones). Constant weight codes possess a rich combinatorial
structure related to the properties of the Johnson graph \cite{del73};
however, no similar theory has arisen for the spherical case.
We provide an explanation of this by showing that the asymptotic
problem of constructing spherical cap codes is equivalent to the
analogous problem for codes on the entire sphere.

Section \ref{sect:hemi} is devoted to a particular case of cap codes, namely
codes in hemispheres. We derive an upper bound on the size of such codes
and use it to derive estimates of a parameter closely related to $k(n),$
the so-called one-sided kissing number $B(n).$
We derive estimates of $B(n)$ for $n=5,6,7,8$ and make a conjecture
about the exact values for some of these cases.

In Section \ref{sect:omega} we use the method of Section \ref{sect:hemi}
to derive another upper bound on spherical cap codes that often improves
the result of Section \ref{sect:bound}. 
We show in examples that the bounds derived in this paper are
sometimes better than the previously known results.
In Section \ref{sect:lp} we present a general approach to bounding the
size of codes in spherical caps that combines several features of the
methods introduced earlier in the paper. We conclude with a brief 
discussion of applications of the bounds on cap codes.

\section{Notation and preliminaries}\label{sect:notation}
Let $S^{n-1}$ be a unit sphere in $n$ dimensions and let 
$e_n=(0,\dots,0,1)$ be the ``North pole.'' Let $0\le \psi\le\phi
\le 90^\circ.$  
Denote by 
   $$
       Z(n,[\psi,\phi])=\{x\in S^{n-1}|\,
                \cos\phi\le\langle x,e_n\rangle\le\cos\psi\}
   $$
a strip cut on the sphere by two planes perpendicular to the vector 
$e_n.$ In particular, $\Cap(n,\phi)=Z(n,[0,\phi])$ is a spherical cap 
with angular radius $\phi$ drawn about $e_n$. 
A finite subset $C\subset Z$ is called a code. Below by $\dist(\cdot,\cdot)$
we denote the angular distance between two points on the sphere.
If a code $C$ has minimum angular 
separation $\theta,$ i.e., satisfies $\dist(x_1,x_2)\ge \theta$ 
for any two distinct points $x_1,x_2\in C,$ we 
call it a $\theta$-{\em code}.
Let 
   $$
       A(n,\theta,[\psi,\phi])=\max_{\substack
     {C \subset Z(n,[\psi,\phi])\\C \text{ a $\theta$-code}}}
              |C|
   $$
be the maximum size of a $\theta$-code in the strip $Z.$ For spherical
caps we will write $A(n,\theta,\phi)$
instead of $A(n,\theta,[0,\phi])$ and use a separate notation 
$B(n,\theta):=A(n,\theta,\nicefrac\pi2)$ for codes in the hemisphere
$S_+:=\Cap(n,\nicefrac\pi2).$ 
In the case of $\phi=\pi$
(the entire sphere) we will call $C$ a spherical $\theta$-code 
and use the notation $A(n,\theta)$ to denote its maximum possible size.

\remove{
Let $C=C_+\cup C_-$ where $C_+$ is a $\theta$-code in the 
hemisphere $S_+$ and $C_-$ is the code obtained from $C_+$
by reflecting its points with respect to the equator plane $x_n=0.$
If $C$ is also a $\theta$-code, we say that $C_+$ can be extended
to an antipodal code. Denote by $B^\ast(n,\theta)$ the maximum size
of a code in $S_+$ that affords such extension.}

The quantity $k(n)=A(n,\nicefrac{\pi}3)$ is equal to the number
of nonoverlapping unit spheres that can touch the sphere $S^{n-1}$ and is 
called the {\em kissing number} in dimension $n$. 
The problem of finding or bounding $k(n)$ has a rich
history \cite{boy94,mus05b,odl79,pfe05}.

Likewise, the quantity
$B(n)=B(n,\nicefrac\pi3)$ is called the {\em one-sided}
kissing number. 
The one-sided kissing number problem was considered recently in
\cite{bez04,bez05}.  $B(n)$ has the following geometric
meaning. Let $H$ be a closed half-space of the $n$-dimensional Euclidean
space. Suppose that $S$ is a unit sphere in $H$ that touches the
supporting hyperplane of $H$. The one-sided kissing number $B(n)$ is
the maximal number of unit nonoverlapping spheres in $H$ that can
touch $S$.

The function $A(n,\theta)$ has received considerable attention
in the literature. Therefore, one possible avenue of studying
spherical cap codes is to map them on the sphere or hemisphere
and relate them to spherical codes.
In this paper, we rely on a
number of mappings between spheres, spherical caps, and spherical strips to
estimate the maximum size of a code in a spherical cap.
Some of them have been used earlier in the literature while the others
have not been emphasized in the context of estimating the code size.
The main problem addressed here is to design the mappings so that
the distance between the images of two points in the domain 
can be bounded in terms of the original distance.
One often-used map is the orthogonal projection $\Pi_n$
which sends the point $x\in S^{n-1}$ along its meridian to the equator of 
the sphere, i.e., the set of points on $S^{n-1}$ with $x_n=0$ 
($\Pi_n$ is defined on $S^{n-1}$ without the North
and South poles). Below we use the notation $S^{n-2}$ to refer to
the equator of the sphere $S^{n-1}.$

Throughout this paper we use the function $\omega(\theta,\alpha,\beta)$
defined by
   $$
     \cos\omega(\theta,\alpha,\beta)=\frac{\cos\theta-\cos\alpha\cos\beta}
{\sin\alpha\sin\beta} .
   $$
In the case of $\alpha=\beta$ we write $\omega(\theta,\alpha)$ instead
of $\omega(\theta,\alpha,\alpha).$
This function describes the change of the distance between two points 
on $S^{n-1}$ which are $\alpha$ and $\beta$ away from $e_n$ and 
$\theta$ away from each other under the action of $\Pi_n.$

\section{Spherical strip (cap) codes and spherical codes}\label{sect:ls}

Several estimates on the size of spherical cap codes have previously
appeared in the literature. They connect 
the maximum size of codes in a spherical cap,
and more generally, in a spherical strip and on the entire sphere.

\subsection{Spherical cap codes and spherical codes} 
Let $m(n,d)$ be the maximum number of points in a unit ball in $\reals^n$
that lie at Euclidean distance $d$ or more apart.
Bounds on $A(n,\theta,\phi)$ are given in the following theorems.

\begin{theorem} \label{thm:cs}
{\rm (Sidelnikov \cite{sid74}, Levenshtein \cite{lev75})} 
   $$
    A(n,\theta,\phi) \le m\Big(n,\frac{2\sin(\theta/2)}{\sin\phi}\Big).
   $$
\end{theorem}

The proof is based on a mapping $\delta:\Cap(n,\phi)\to\reals^n$
that transforms the cap to the unit ball in $\reals^n$
according to the following rule:
  $$
   \delta(x)=\frac1{\sin\phi} (x-e_n\, \cos\phi).
  $$

\remove{The mapping $Q$ projects the vector $(1/\sin \phi)\delta(x)$
on the sphere  in the next dimension 
(i.e., sends the equator plane of $S^n$ to the hemisphere along the direction
$e_{n+1}$ orthogonal to it).}

\begin{theorem}\label{thm:lv} {\rm (Levenshtein \cite{lev83a}).} 
   $$
     m\Big(n-1, \frac{2\sin\theta/2}{\sin\phi\cos\phi}\Big)
       \le A(n,\theta,\phi)\le 
        m(n-1,2\sin\nicefrac\theta2\cot\phi).
   $$
\end{theorem}
The proof is based on a mapping that projects the cap centrally on the
tangent hyperplane to the sphere $S^{n-1}$ at the point $e_n.$

Bounds of these two theorems 
are useful in asymptotics (both as $n\to\infty$ and as
$\theta\to 0$) for estimating the size of spherical codes \cite{kab78}
and the packing density in $\reals^n$ \cite{sid74,lev83a}.
Their use for finite $n$ is based on the obvious inequality
   $
    m(n-1,2d)\le A(n,2\arcsin d)
  $
and leads to the estimates
  \begin{align}
    A(n,\theta,\phi)&\le A(n+1,2\arcsin({\sin\nicefrac\theta2}\,{\cosec\phi}))
         =A(n+1,\omega(\theta,\phi))
            \label{eq:ls1}\\
    A(n,\theta,\phi)&\le A(n,2\arcsin({\sin\nicefrac\theta2}\,{\cot\phi})).
          \label{eq:ls2}
  \end{align}
For large $n$ upper bound (\ref{eq:ls1}) is uniformly
better than bound (\ref{eq:ls2}) because 
$\nicefrac1{\sin\phi}>{\cot\phi}.$ 
For finite $n$ bound (\ref{eq:ls2})
is stronger than (\ref{eq:ls1}) for small $\phi$ and is weaker than it
otherwise.

\subsection{Spherical strip codes}

\begin{figure}[H]
\begin{center}\includegraphics[width=2.5in]{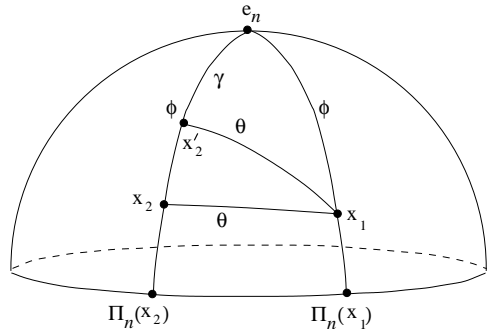}
\end{center}
\caption{Definition of the angle $\gamma$: the triangle $x_1x_2x_2'$ is 
isosceles.}\label{fig:iso}
\end{figure}

In this subsection we discuss the action of the projection $\Pi_n$
on the code in a spherical strip $Z(n,[\psi,\phi]).$
Given a $\theta$-code 
$C\subset Z$ we would like to know what happens to its distance
upon applying the mapping $\Pi_n$ to it.
Given two points $x_1,x_2\in Z$ the main issue is to establish 
how the distance between their images depends on their relative location 
in the strip.
Let $\dist(x_1,x_2)=\theta$ and let the angle
$\gamma\ne \phi, 0\le\gamma\le 90^\circ$ be defined by the equation
    $$
      \cos\omega(\theta,\phi)=\cos\omega(\theta,\phi,\gamma).
    $$
Geometrically, the angle $\gamma$ is defined as follows. Consider 
two
points $x_1,x_2$ that are $\theta$ away from each other and
lie on the boundary of the cap $\Cap(n,\phi)$ (i.e., the angle between
each of them and $e_n$ is $\phi$). The distance
between their images under $\Pi_n$ equals $\omega(\theta,\phi).$
Consider the point $x_2'$ that satisfies $\Pi_n(x_2)=\Pi_n(x_2')$
and $\langle x_2,x_2'\rangle=\cos\theta,$ then $\gamma=\dist(x_2,x_2')$
(see Figure \ref{fig:iso}). Formally, $\gamma$ is the angle given by
  $$
    \sin\gamma=\frac{\sin\phi(\cos^2\theta-\cos^2\phi)}{\cos^2\theta
      +\cos^2\phi(1-2\cos\theta)}.
  $$
Under the mapping $\Pi_n$ the code in a strip is transformed to a
spherical code in $S^{n-2}$. This transformation yields nontrivial
results only in the case of $\theta>\psi-\phi$ (otherwise the code $C$
can contain points that project identically on the equator, so the
distance of the image code is zero).
In this case, the distance in the image code is minimized for a
pair of points on the ``lower'' boundary of the strip if 
$\gamma<\psi$ and for a pair of points one of which is on the
lower and the other on the upper boundary, otherwise.

More precisely, we have the following theorem.
\begin{theorem} \label{thm:avz} {\rm (Agrell, Vardy and Zeger~\cite{agr00})}
Let $0\le \psi\le\phi\le \nicefrac\pi2$ and $\theta<2\phi.$

$(a)$ Let $\theta>\phi-\psi.$
If $\gamma<\psi$ then
     $$
A(n,\theta,[\psi,\phi])=A(n-1,\omega(\theta,\phi)).$$
If $\gamma>\psi$ then
   \begin{gather}\label{eq:caseb}
      A(n,\theta,[\psi,\phi])\le
         \min\Big\{A(n,\theta,[\psi,\gamma]) + A(n-1,\omega(\theta,\phi)),
          A(n-1,\omega(\theta,\phi,\psi))\Big\}
   \end{gather}

$(b)$ Let $\theta\le \phi-\psi.$ Then
   $$
    A(n,\theta,[\psi,\phi])\le A(n,\theta,[\psi,\gamma]) + 
      A(n-1,\omega(\theta,\phi))            \quad(\gamma>\psi)
   $$
If $\theta>2\phi$ then $A(n,\theta,[\psi,\phi])=1.$
\end{theorem}

By taking $\psi=0$ this theorem implies the following corollary for 
spherical cap codes.
\begin{corollary}\label{cor:agr00} {\rm \cite{agr00}}
\begin{equation}\label{eq:capsavz}
    A(n,\theta,\phi)\begin{cases}
      \le A(n,\theta,\gamma)+ 
            A(n-1,\omega(\theta,\phi))
         \quad &{\text{ {\rm if} }} 0<\theta\le \phi\\
      = A(n-1,\omega(\theta,\phi))
         &{\text{ {\rm if} }} \phi<\theta\le 2\phi\\
      =1 &{\text{ {\rm if} }} \theta > 2\phi.
   \end{cases}
\end{equation}
\end{corollary}        

\subsection{Lower bounds} A general lower (existence) bound on the size of a
$\theta$-code in the cap $\Cap(n,\phi)$ can be obtained by the standard
greedy argument. It follows that there exist codes of size 
   $$
    M\ge \Omega_n(\phi)/\Omega_n(\theta),
   $$
where by $\Omega_n(\beta)=\frac{2\pi^{(n-1)/2}}{\Gamma((n-1)/2)}\int_0^\beta
\sin^{n-2}\tau d\tau$ we denote the area of the spherical cap 
$\Cap(n,\beta)$ on the sphere $S^{n-1}$.

\section{An extension of Delsarte's method}\label{sect:extension}

In this section we explain a way to use bounds on spherical cap codes
in order to extend the well-known Delsarte method for bounding the size
of spherical codes.

The original Delsarte 
(linear programming) bound as applied to spherical codes \cite{del77b,kab78}
has the following form. Let $\{G_k^{(n)}(x)\}_{k=0,1,\dots}$ denote the 
family
of Gegenbauer polynomials, i.e., polynomials orthogonal on $[-1,1]$
with weight $(1-x)^{(n-3)/2}$ and satisfying the normalization 
condition $G_k^{(n)}(1)=1.$ Suppose that a real function $f$ is
a nonnegative linear combination of Gegenbauer polynomials $G_k^{(n)}(t)$,
i.e.,
       $$
       f(t)=\sum\limits_k {f_kG_k^{(n)}(t)},\; \mbox{ where }\; f_k\ge 0.
       $$  
If $f(t)\le 0$ for all $t\in [-1,\cos\theta]$ and $f_0>0,$ 
then $A(n,\theta)\le f(1)/f_0.$

Next we consider an extension of this method to spherical caps.
Let $Y=\{y_1,\ldots,y_m\}$ be a $\theta$-code in the spherical cap
$\Cap(n,\phi)$ with center $e_n$ and let $\cY$ be the set of all such codes.
Of course, $m\le A(n,\theta,\phi).$ Let $e_n^\ast=-e_n$,
let $f(t)$ be a real function on the interval $[-1,1]$,
         $$H_f(Y)=H_f(y_1,\ldots,y_m):=f(1)+f(\langle
        e_n^\ast,y_1\rangle)+\ldots+f(\langle e_n^*,y_m\rangle),
           $$
           $$h_m(n,\theta,\phi,f):=\max\limits_{Y\in \cY}\{H_f(Y)\},\quad
            h_{\max}(n,\theta,\phi,f):=\max\limits_{m\le
                A(n,\theta,\phi)}{\{h_m(n,\theta,\phi,f)\}}.
           $$

\begin{theorem} \label{thm:extension}
Suppose that $f$ is a nonnegative linear combination of 
Gegenbauer polynomials $G_k^{(n)}(t)$, i.e.,
       $$
       f(t)=\sum\limits_k {f_kG_k^{(n)}(t)},\; \mbox{ where }\; f_k\ge 0.
       $$
If $f(t)\le 0$ for all $t\in [-\cos{\phi},\cos{\theta}]$ and $f_0>0$,  then
        $$
          A(n,\theta)\le \frac{h_{\max}(n,\theta,\phi,f)}{f_0}.
       $$
\end{theorem}

\begin{proof} Let $C=\{x_1,\ldots,x_M\}$ be a 
$\theta$-code in  ${S}^{n-1}$.
It is well known \cite{del77b,kab78} that
$$\sum\limits_{i=1}^M\sum\limits_{j=1}^M {{G_k^{(n)}(t_{i,j})}} \ge 0, 
\quad t_{i,j}:=\langle x_i, x_j\rangle =\cos{(\dist(x_i,x_j))}.$$

Using this we obtain
$$S_f(C):=\sum\limits_{i=1}^M\sum\limits_{j=1}^M{f(t_{i,j})}=
\sum\limits_{k,i,j}{f_kG_k^{(n)}(t_{i,j})}\ge  
\sum\limits_{i,j}{f_0G_0^{(n)}(t_{i,j})} =   f_0M^2.$$

Let
      $$
       J(i):=\{j:f(\langle x_i,x_j\rangle)>0, \; j\neq i\}, 
\quad C(i)=\{x_j\in C: j\in J(i)\}, m_i=|C(i)|.
         $$

Note that $j\in J(i)$ only if $x_j$ belongs to the $\Cap(n,\phi)$ 
with the center at $-x_i$.
Then
          $$
           S_i(C):=\sum\limits_{j=1}^M {f(\langle x_i,x_j\rangle)}
          \le f(1)+\sum\limits_{j\in J(i)} 
               {f(\langle x_i,x_j\rangle)}
            =H_f(C(i))\le h_{m_i}(n,\theta,\phi,f).
           $$ 
Therefore, 
       $$
        f_0M^2\le S_f(C)=\sum\limits_{i=1}\limits^M S_i(C)\le Mh_{\max},
        $$
i.e. $M\le h_{\max}/f_0$ as required.
\end{proof}

Note that $h_0=f(1).$ If $f(t)\le 0$ for all 
$t\in [-1,\cos{\theta}]$, then all
$m_i=0,$ i.e. $h_{\max}=h_0=f(1)$ and $M\le f(1)/f_0,$
so this theorem includes the Delsarte bound as a particular case.

For given $n, \theta, \phi, f$ and $m$ the value
$h_m(n,\theta,\phi,f)$ is the solution of the following optimization
problem on $S^{n-1}$:
           $$
h_m(n,\theta,\phi,f)=f(1)+f(-\cos{\phi_1})+\ldots+f(-\cos{\phi_m})\;
\to \; \max
          $$ 
subject to the constraints
       $$
        \phi_i:=\dist(e_n^*,y_i)\le\phi,\; 1\le i\le m; \quad
        \dist(y_i,y_j)\ge \theta,\; i\ne j.
        $$ 
The dimension of this problem is $m(n-1)\le(n-1)A(n,\theta,\phi)$.  
For relatively small $n$ and $A(n,\theta,\phi)$ optimization can
be carried out numerically.  
Moreover, if in addition to the above restrictions the function $f(t)$ 
is monotone decreasing for 
$t\in [-1,-\cos{\phi}]$ then in some cases the dimension of this problem 
can be reduced to $n$ (see the details in \cite{mus04,mus05b}).
 Suitable polynomials
$f$ can be found by linear programming (see an algorithm in the
Appendix to \cite{mus05b}).

\section{The case of large angles}\label{sect:large}
In this section we consider $\theta$-codes in a spherical
 cap $\Cap(n,\phi)$ with large values of $\theta$. More precisely
let us assume that $\theta>\phi.$
Clearly, if $\theta>2\phi$, then no more than one point can lie in
$\Cap(n,\phi)$, i.e.  $A(n,\theta,\phi)=1$. Now we consider the case
$2\phi\ge\theta>\phi.$ Recall that $\phi\le \pi/2.$

\remove {It was proved (see [ ]) that points in such codes form a 
spherical convex polytope, i.e.,
 these points are vertices of its convex hull.}

\begin{lemma}\label{lemma:60} Suppose $\; 2\phi\ge\theta>\phi>0$, then $\; 
\omega(\theta,\phi)>\pi/3$. 
\end{lemma}

\begin{proof} Let $z:=\cos{\theta},\; t:=\cos{\phi}$. Then
       $$\cos{\omega(\theta,\phi)}=
\frac{\cos{\theta}-\cos^2{\phi}}{\sin^2{\phi}}=\frac{z-t^2}{1-t^2}\le
\frac{z-z^2}{1-z^2}=\frac{z}{1+z}<\frac{1}{2}.$$
Thus, $\omega(\theta,\phi)>\pi/3$.
\end{proof}

As stated in Theorem \ref{thm:avz}, for $2\phi\ge\theta>\phi$
the problem of finding $A(n,\theta,\phi)$ is equivalent to bounding
the size of spherical codes.
Since the proof in \cite{agr00} is not isolated into a separate argument
we include it here for completeness.

\begin {theorem}\label{thm:exact} 
{\rm (Agrell et al.~\cite{agr00}, Musin~\cite{mus05b})} 
If $\; 2\phi\ge\theta>\phi$, then
$$A(n,\theta,\phi)=A(n-1,\omega(\theta,\phi)).$$ 
\end {theorem}
\begin{proof} 
First let us prove the lower bound.
\begin{lemma} \label{lemma:lower}
Let $0\le\psi<\phi\le \pi/2,\;\theta\le 2\phi.$ Then
    $$
      A(n,\theta,\phi)\ge A(n,\theta,[\psi,\phi])\ge 
          A(n-1,\omega(\theta,\phi)).
    $$
\end{lemma}
\begin{proof}
The first inequality is obvious. To prove the second one let
consider the strip $Z(n,[\psi,\phi])\subset S^{n-1}$ and let $\Sigma$ be its 
``lower'' boundary. The projection $\Pi_n$ is a one-to-one map from 
$\Sigma$ to the unit sphere $S^{n-2}$ 
(the equatorial sphere of $S^{n-1}$). Now consider a code $C'\subset S^{n-2}$
and the code $C\subset \Sigma$ that corresponds to $C'$ under this map.
If the distance of $C'$ equals $\omega(\theta,\phi),$ then the
distance of $C$ is $\theta$ (the function $\omega(\theta,\phi)$ is
monotone). Since $|C'|=|C|\le A(n,\theta,[\psi,\phi]),$
this proves the needed inequality.
\remove{
$C\subset S^{n-2}$ be a spherical code with distance $\omega(\theta,\phi).$
Let $\Sigma$ be 
$C\subset \Cap(\theta,\phi)$ be a $\theta$-code.
 Let $\Sigma$ be the ``lower'' boundary  of the strip, i.e., 
a sphere in $n-1$ dimensions of radius $\sin\phi$, and let $C'$ be a 
spherical $\theta$-code on $\Sigma.$ It is easy to see that under the 
map which sends 
$\Sigma$ to the unit sphere $S^{n-2}$ by increasing the distances
by $1/\sin\phi,$ the code $C'$ becomes an $\omega(\theta,\phi)$-code.
Since $|C|\ge |C'|,$ this proves the lower bound.
To prove the second one, 
let $x_i,x_j\in C$ be two
code vectors on the ``lower'' boundary $\Sigma$ of the strip (i.e., satisfying
$\dist(x_1,e_n)=\dist(x_2,e_n)=\phi$ ) and suppose that
$\dist(x_1,x_2)=\theta.$  By (\ref{eq:cos})
   \begin{equation}\label{eq:sc}
      \cos\theta=\cos^2\phi+\sin^2\phi \cos\gamma_{i,j},
   \end{equation}
so $\gamma_{i,j}=\omega(\theta,\phi).$
Note that for any fixed $\theta$ the function $\omega(\theta,\phi)$
is monotone (increasing) on $\phi.$
Therefore, under the inverse projection $S^{n-2}\to\Sigma$ the code
$C'\in S^{n-2}$ with distance $\omega(\theta,\phi)$ will
be transformed to a $\theta$-code $\Sigma \Cap(n,\phi).$
Thus, 
$$A(n,\theta,[\psi,\phi])\ge A(n-1,\omega(\theta,\phi)).$$}
\end{proof}

Now let $C=\{x_1,\ldots,x_m\}$ be a $\theta$-code in  
$\Cap(n,\phi)$. Then
$$\; \theta_{i,j}:=\dist(x_i,x_j)\ge \theta \;\; \mbox{ for } \; i\neq j.$$
Denote by $\phi_i$ the angular distance between $e_n$ and $x_i$, where $e_n$ 
is the center of  $\Cap(n,\phi)$. Note that $\phi_i\le \phi.$
 
Let $X=\Pi_n(C)$ be the image of $C$ under the projection on the equator 
$S^{n-2}$ of the sphere from its North pole $e_n.$ 
Denote by  $\gamma_{i,j}=\dist(\Pi_n(x_i),\Pi_n(x_j))$ be the distance
between the images of $x_i$ and $x_j$ under the projection.
\remove{\begin{figure}[H]\begin{center}
\includegraphics[width=3in]{cap2.eps}\end{center}
\caption{}\end{figure}}
Recall the law of cosines for a spherical triangle. Suppose the two
sides are $a,b$ and the angle between them is $\psi$, then the
third side $c$ satisfies
  \begin{equation}\label{eq:cos}
    \cos c=\cos a\cos b+\cos\psi\sin a\sin b.
    \end{equation}
From this and the inequality 
$\cos{\theta_{i,j}}\le \cos{\theta},$ we get
$$\cos{\gamma_{i,j}}=\frac{\cos{\theta_{i,j}}-
\cos{\phi_i}\cos{\phi_j}}{\sin{\phi_i}\sin{\phi_j}}
\le 
\frac{\cos{\theta}-\cos{\phi_i}\cos{\phi_j}}{\sin{\phi_i}\sin{\phi_j}}
$$ 
$$\mbox{Let }\quad Q(\alpha,\beta)=\frac{\cos{\theta}-\cos{\alpha}
\cos{\beta}}{\sin{\alpha}\sin{\beta}}, \; \; \mbox{ then } \; \;
\frac{\partial Q(\alpha,\beta)}{\partial\alpha}=\frac{\cos{\beta}-
\cos{\theta}\cos{\alpha}}{\sin^2{\alpha}\sin{\beta}}.$$ 
From this it follows that if 
$\; 0<\alpha, \beta\le \phi$ then 
$ \cos{\beta}\ge \cos{\theta}$ (because $\theta\ge \phi$); therefore 
${\partial Q(\alpha,\beta)}/{\partial\alpha}\ge 0$, i.e., $Q(\alpha,\beta)$ 
is a monotone increasing function in $\alpha$. We have 
$Q(\alpha,\beta)\le Q(\phi,\beta)=Q(\beta,\phi)\le Q(\phi,\phi).$ Therefore,
$$\cos{\gamma_{i,j}}\le \frac{\cos{\theta}-\cos{\phi_i}\cos{\phi_j}}
{\sin{\phi_i}\sin{\phi_j}} \le
\frac{\cos{\theta}-\cos^2{\phi}}{\sin^2{\phi}}=\cos{\omega(\theta,\phi)}.$$
Thus $X$ is an $\omega(\theta,\phi)$-code on the $(n-2)$-sphere.  
That yields $$A(n,\theta,\phi)\le A(n-1,\omega(\theta,\phi)).$$
\end{proof}

It is proved in \cite{mus05c,mus04} that in the case covered by this theorem,
points in an extremal configuration are vertices of a convex polyhedron,
and lie on the boundary of the cap. This implies that if $\theta=\phi$ then
the code can be augmented by the point $e_n$ without reducing its
distance, so $A(n,\theta,\theta)=A(n,\omega(\theta,\theta))+1.$

Denote by $\varphi_n(M)$ the largest angular distance in a spherical
code on $S^{n-1}$ that contains $M$ points.
Recall that $k(n)$ denotes the kissing number in $n$ dimensions.

\begin{corollary} Suppose that $\theta>\phi$, then 
$$A(n,\theta,\phi)\le k(n-1).$$ Moreover, if $\; \varphi_{n-1}(K)\le\pi/3$,
then $$A(n,\theta,\phi) < K.$$
\end{corollary}

\begin{proof} By Lemma \ref{lemma:60} we can write
$\; \omega(\theta,\phi)=\pi/3+\varepsilon, \; \varepsilon>0$. 
Then the theorem yields $$A(n,\theta,\phi) = 
A(n-1,\omega(\theta,\phi))=A(n-1,\pi/3+\varepsilon)\le A(n-1,\pi/3)=k(n-1).$$
If $\; \varphi_{n-1}(K)\le\pi/3$, then $A(n-1,\pi/3+\varepsilon)<K$.
\end{proof}

Using Theorem \ref{thm:exact} together with this
corollary we can find the exact 
value $A(3,\theta,\phi)$ and $A(4,\theta,\phi)$ for $\theta>\phi$.
We can also find $A(n,\theta,\phi)$ for all $n$ if $\cos{\theta}<\cos^2{\phi}.$

\begin{enumerate}
\item
  Let $n=3$. Note that $k(2)=6, \; \varphi_2(6)=\pi/3,$ and
$\; \varphi_2(M)=2\pi/M.$ Then
$$A(3,\theta,\phi)=\lfloor 2\pi/\omega(\theta,\phi)\rfloor\le 5.$$

\item Let {$n=4$.} In three dimensions the best codes and the 
values $\varphi_3(M)$ are presently known for $M\leqslant 12$   and $M=24$ 
(see \cite{dan86,fej53,sch51}).
It follows from Fejes T\'oth's bound \cite{fej53} that 
       $$\varphi_3(2)=180^\circ, \quad \varphi_3(3)=120^\circ, 
         \quad \varphi_3(4)=\arccos(-1/3)\approx 109.47^\circ,$$
            $$\varphi_3(6)=90^\circ, 
        \quad \varphi_3(12)=\arccos{(1/\sqrt{5})}\approx 63.435^\circ.$$

Sch\"utte and van der Waerden \cite{sch51} proved that 
$$\varphi_3(5)=\varphi_3(6)=90^\circ, \; \quad\varphi_3(7)
\approx 77.87^\circ \; (\cos{\varphi_3(7)}=\cot{40^\circ}\cot{80^\circ}),$$ 
$$\varphi_3(8)=\arccos{\frac{\sqrt{8}-1}{7}}\approx  74.86^\circ,\qquad 
\varphi_3(9)=\arccos{\frac{1}{3}}\approx  70.53^\circ.$$

The cases $M=10, 11$ were considered by Danzer \cite{dan86}:
$$\varphi_3(10)\approx 66.15^\circ, \qquad \varphi_3(11)=\varphi_3(12).$$

Since $k(3)=12$ \cite{sch53}, we have $\; A(4,\theta,\phi)\le 12.\; $ Thus
$$A(4,\theta,\phi)=\max\limits_{M\le 12}{\{M: \; 
  \varphi_3(M)\ge \omega(\theta,\phi)\}}.$$

\item
Let $\cos{\theta}<\cos^2{\phi}.$ In this case we have 
$\omega(\theta,\phi)>90^\circ.$
It is well known \cite{ran55} that for all dimensions 
$$\varphi_n(M)=\arccos{\left(\frac{-1}{M-1}\right)}, \quad 2\le M \le n+1;$$
and
$$\varphi_n(n+2)=\ldots=\varphi_n(2n-1)=\varphi_n(2n)=90^\circ.$$
Therefore, for $\; \arccos(\cos^2{\phi})<\theta\le\ 2\phi$, we obtain:
$$A(n,\theta,\phi)=\max\limits_{M\le n}{\{M:\; 
\arccos{\frac{-1}{M-1}}\ge \omega(\theta,\phi)\}}.$$
\end{enumerate}

The results obtained can be applied for the kissing number problem as follows.

For $n=3$ let us consider the following polynomial $f$:
        $$
         f(t) = \frac{2431}{80}t^9 - \frac{1287}{20}t^7 +
\frac{18333}{400}t^5 + \frac{343}{40}t^4 - \frac{83}{10}t^3 -
\frac{213}{100}t^2+\frac{t}{10} - \frac{1}{200}.
      $$ 
This polynomial 
satisfies the assumptions of Theorem \ref{thm:extension}
 with $\theta=60^\circ,$ 
$\phi\approx 53.794^\circ$ $ (f(-\cos{\phi})=0),$ and $f_0=1.$ In this
case $A(3,\theta,\phi)=4.$ Since $h_{\max}<13$ (see a proof in
\cite{mus05a}) we have $k(3)=A(3,\pi/3)<13$, i.e. $k(3)=12.$

In the case of $n=4,$ Theorem \ref{thm:extension} can be applied with
        $$
        f(t) = 53.76t^9 - 107.52t^7 + 70.56t^5 + 16.384t^4 - 
          9.832t^3 - 4.128t^2 -
         0.434t - 0.016.
           $$
Here $\theta=60^\circ, \; \phi\approx 52.559^\circ, f_0=1$, and 
$A(4,\theta,\phi)=6.$ It was proved \cite{mus05b} that $h_{\max}<25$. 
Since $k(4)\ge 24$ this yields $k(4)=24.$

Recently, Pfender \cite{pfe05} considered the case 
$\cos{\theta}<\cos^2{\phi}.$ 
He found some improvements for upper bounds on $k(n)$ for dimensions
$n=9, 10, 16, 17, 25, 26.$

\section{Stretching transformation}\label{sect:bound}

In this section we will prove the following bound on spherical cap
codes which relates the maximum size of such a code to the size of
codes in a hemisphere.
 \begin{theorem}\label{thm:hemi} Let $\nicefrac\theta 2< \phi\le \nicefrac
\pi2.$ Then
    \begin{equation}\label{eq:hemic}
          A(n,\theta,\phi) \le B(n,\omega(\theta,\phi)).
    \end{equation}
 \end{theorem}

{\em Remark.} 
This theorem improves upon the bound (\ref{eq:ls1})
by reducing the dimension on the right-hand side by one. 
It also extends
the applicability of the bound in Theorem \ref{thm:exact} to the range
of angles $\theta/2\le\phi\le\theta,$ although in this range we cannot
claim the exact equality anymore. On the other hand, by Theorem \ref{thm:hemi}
it is sufficient to estimate the number of code points in the hemisphere
as opposed to the entire sphere.
 
Note also that for any $\theta,\phi,$
  $$
    \cos\omega(\theta,\phi)-\cos\theta=\tan^2\phi(\cos\theta-1)<0
  $$
therefore, this theorem is stronger in the entire range of angles
than the trivial bound 
$A(n,\theta,\phi)\le A(n,\theta)$. Finally, 
the angle $\omega(\theta,\phi)$ ranges between $\pi$ and
$\theta$ as $\phi$ grows from $\nicefrac\theta 2$ to $\nicefrac \pi 2$ and
is a monotone decreasing function of $\theta.$ 

\medskip
A proof of Theorem \ref{thm:hemi} 
will follow from the following result which 
describes the effect on the distance of spherical cap codes
of a ``stretching map'' of spherical caps. 

\begin{theorem}\label{thm:recursion} Let $0< \phi\le \pi/2.$ Then
for any $s\ge 1$
   $$
     A(n,\theta,\phi)\le A(n,\theta',s\phi),
   $$
where 
  $$
   \cos \theta'=\cos^2 s\phi + \frac{\sin^2 s\phi}
          {\sin^2 \phi}(\cos\theta-\cos^2\phi).
  $$
\end{theorem}
{\em Proof:}
Let $T_s, s\ge 1$ be a map
on $\Cap(n,\phi)$ defined as follows: for a point 
$x\in \Cap(n,\phi)$ with $\dist(x,e_n)=\alpha$ its image 
$\bfy=T_s(x)$
is a point that satisfies $\dist(\bfy,e_n)=s\alpha$ 
and lies on the meridian passing through $x.$
Thus, $T_s(\Cap(n,\phi))=\Cap(n,s\phi),$ and we assume that 
$s\phi\le \nicefrac\pi2.$

The proof of Theorem \ref{thm:recursion} relies upon
the next two lemmas.
\begin{lemma} Let $x_1,x_2\in \Cap(n,\phi)$ with 
$\dist(x_1,e_n)=u,\dist(x_2,e_n)=v,\linebreak[4]\dist(x_1,x_2)=\theta.$
The distance $\dist(T_s(x_1),T_s(x_2))$ reaches its minimum when $u=v.$
\end{lemma}
\begin{proof} Figure \ref{fig:lemma} shows the relative location on the sphere
of $x_1,x_2$ and their images $\bfy_1=T_s(x_1),\bfy_2=T_s(x_2).$
\begin{figure}[H]\begin{center}\includegraphics[width=2.5in]{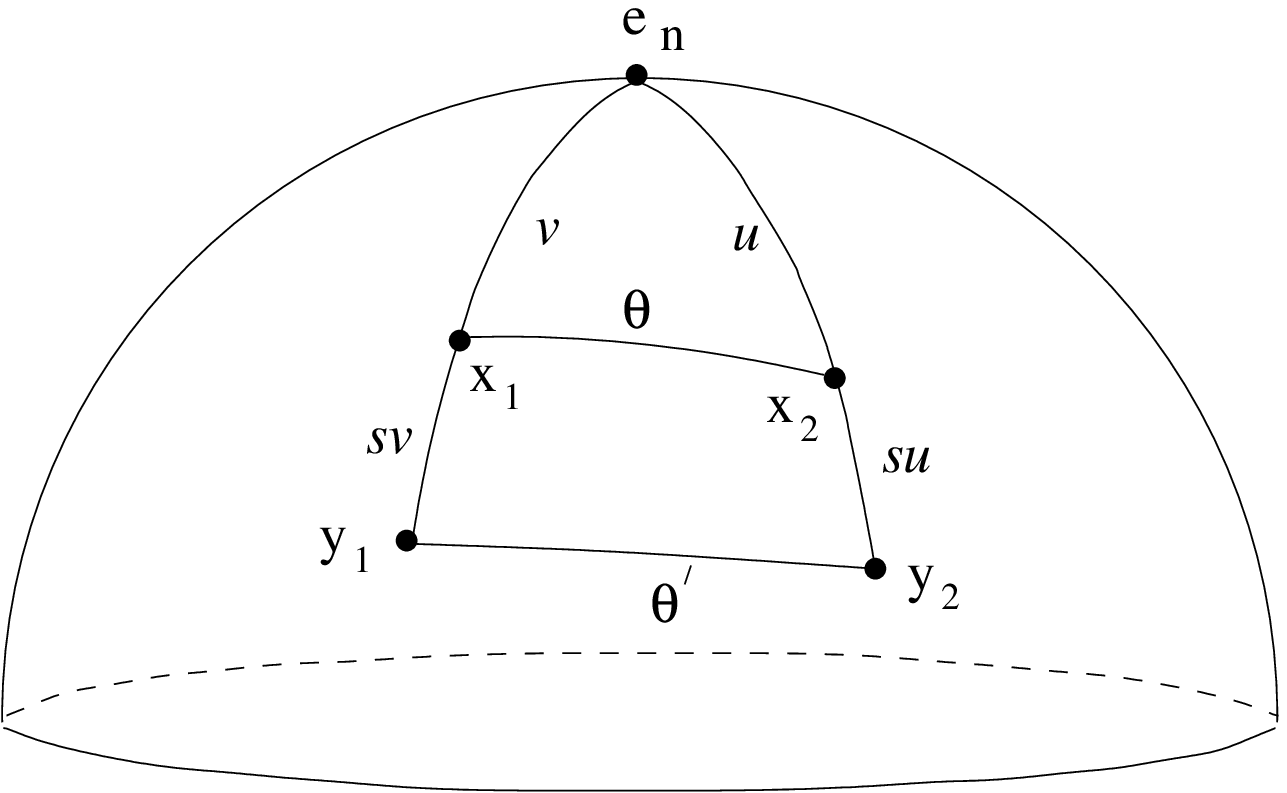}
\end{center}
\caption{}\label{fig:lemma}
\end{figure}
\noindent Using (\ref{eq:cos}) we find that $\cos\theta'=F(u,v,\theta),$
where
  $$
   F(u,v,\theta)=\cos su\cos sv+\rho(u)\rho(v)(\cos\theta-\cos u \cos v),
  $$
$\rho(t)=\sin(st)/\sin t.$ For definiteness assume that $v\le u.$
We need to prove that 
$F(u,v,\theta)\le F(u,u,\theta)$ where $0\le u-v\le\theta$ and
$su\le \pi/2.$

{\em Fact 1.} (i) The function $\rho(t)$ is monotone decreasing for
$t\in (0,\pi/2s).$ Indeed 
  $$\rho'(t)=\frac{\cos t\cos st}{\sin^2 t}(s\tan t-\tan st)\le 0$$
with the equality only if $t=0$.

(ii) $\rho(t)<s$ (follows from (i) and the equality $\rho(0)=s$).

{\em Fact 2.} The function
   $
     S(\theta)=F(u,v,\theta)-F(u,u,\theta)$
is maximized on $\theta$ for $\theta=u-v.$

{\em Proof:} The coefficient of $\cos \theta$ in $S(\theta)$ equals
$\rho(u)(\rho(v)-\rho(u))\ge 0$. Then the claim follows from the condition
$\nicefrac\pi2\ge\theta\ge u-v.$

{\em Fact 3} (which implies the lemma). $S(\theta)\le S(u-v)\le 0.$

{\em Proof:} The first inequality is proved in Fact 2. Now compute
  $$
   S(u-v)=\cos s(u-v)-\cos^2 su-\frac{\sin^2 su}{\sin ^2u}
               (\cos(u-v)-\cos^2 u)
   $$
  $$
    =\cos s(u-v)-1-\rho^2(u)\cos(u-v)+\rho^2(u).
  $$
The derivative of the last expression on $v$ equals 
$(s \rho(u-v)-\rho^2(u))\sin(u-v)$. Since $\rho(u)\le \rho(u-v)$, we
can write 
  $$
    (s \rho(u-v)-\rho^2(u))\sin(u-v)\ge (s-\rho(u))\rho(u-v)\sin (u-v)\ge 0
  $$
where the last inequality follows from part (ii) of Fact 1 above.
\end{proof}

\begin{lemma} Let $x_1,x_2\in \Cap(n,\phi)$ with 
$\dist(x_1,e_n)=\dist(x_2,e_n)=u,\dist(x_1,x_2)=\theta,$
where $u\le\phi.$
Then the
 distance $\dist(T_s(x_1),T_s(x_2))$ reaches its minimum when $u=
\phi.$
\end{lemma}
\begin{proof}
Since 
   $$
     \cos\theta'=F(u,u,\theta)=1-\rho^2(u)(1-\cos\theta),
   $$
the claim is implied by Fact 1(i) above. \end{proof}

The last two lemmas imply Theorem \ref{thm:recursion}.
Indeed, let $x_1,x_2\in \Cap(n,\theta,\phi)$ be two points
at distance $\theta.$ The lemmas show that in order for the distance
between their images under $T_s$ to reach its minimum the points should
lie on the boundary of the cap. Then the expression for $\theta'$ in
the theorem is implied by an application of the cosine law (\ref{eq:cos}).

Finally, Theorem \ref{thm:hemi} follows from Theorem \ref{thm:recursion}
by taking $s=\pi/2\phi.$ \qed

We conclude this section with two applications of Theorem \ref{thm:hemi}.

\subsection{An upper bound on spherical codes}
Here we establish the following new estimates:
    \begin{equation}\label{eq:be}
       A(n,\theta)\le \frac{\Omega_n}{\Omega_n(\phi)} 
       B(n,\omega(\theta,\phi)).
    \end{equation} 
    \begin{equation}\label{eq:be1}
        A(n,\theta) \le
          \frac{\Omega_n}{\Omega_n(\phi)} A(n-1,\omega(\theta,\phi))
       \quad(\theta>\phi)
    \end{equation}
where $\Omega_n(\phi)$ is the area of the cap of radius $\phi$
and $\Omega_n=\pi^{n/2}/\Gamma(\nicefrac n2+1)$ 
is the ``surface area'' of the unit sphere $S^{n-1}.$
They are implied by the Bassalygo-Elias inequality stated in the next lemma.
\begin{lemma}\label{lemma:be} \cite{sid74,lev75}
   Let $\nicefrac\theta2\le\phi\le\nicefrac\pi2.$ Then
        \begin{equation*}
     A(n,\theta)\le \frac{\Omega_n}{\Omega_n(\phi)}A(n,\theta,\phi).
          \end{equation*}
\end{lemma}
{\em Proof:}~
   Consider a code $C\subset S^{n-1}$ and let $C_\phi(z)$ be the number
of code points in the cap with ``center'' $z$ and radius $\phi.$ 
Note that every cap whose center $z$ is at most $\theta$ away from 
a given code point $x$ will contain this point. Then clearly 
    $$\int_{z\in S^{n-1}}  |C_\phi(z)|dz=\Omega_n(\phi)|C|.$$
Since $C_\phi(z)\le A(n,\theta,\phi),$ we obtain
    $$
      A(n,\theta,\phi)\Omega_n\ge A(n,\theta)\Omega_n(\phi).
    $$
{\hfill\qed}

 Therefore, using Theorems \ref{thm:hemi} and \ref{thm:exact} 
we obtain the bounds (\ref{eq:be}),\,(\ref{eq:be1}).
In particular, 
inequality (\ref{eq:be}) is stronger than bounds on cap codes based
on Lemma \ref{lemma:be} that appeared
in \cite{sid74,lev75,lev83a}.
\remove{Either of these inequalities
leads to the asymptotic improvement of the bound on the rate of spherical
codes with a given distance $\theta$ for low values of $\theta$ obtained
in Kabatyansky and Levenshtein \cite{kab78}. More precisely, bounding
$B(n,\omega(\theta,\phi))$ by $A(n,\omega(\theta,\phi))$ in (\ref{eq:be})
and using the upper bound on spherical codes derived by linear programming,
they derive an improved asymptotic bound on spherical codes for
$\theta\le 63^\circ.$}

\subsection{Large dimensions}
Let $R(C)=\nicefrac 1n \ln |C|$ be the {\em rate} of the code 
$C\subset S^{n-1}.$ Denote by
         $$
      R^+(\theta)=\limsup_{n\to\infty} \frac1n \ln
             A(n,\theta) \quad R^-(\theta)=\
            \liminf_{n\to\infty} \frac1n \ln A(n,\theta).
         $$
Abusing notation, we write below $R(\theta)$ to refer to the common
value of $R^+$ and $R^-$ even though it is not known that the limit exists.
Likewise we write $R(\theta,\phi)$ and $R(\theta,[\psi,\phi])$
to refer to cap and strip codes. In this section we show that the
problem of finding either of the last two quantities 
is equivalent to that of computing
$R(\theta)$. More precisely, we have the following theorem.
\begin{theorem}\label{cor:eq} 
Let $0\le\psi<\phi\le\pi/2.$
Then
      $$
        R(\theta,[\psi,\phi])=R(\theta,\phi)=R(\omega(\theta,\phi)).
      $$
\end{theorem}

Indeed, combining Lemma \ref{lemma:lower}, Theorem \ref{thm:hemi} and
the obvious
inequality $B(n,\alpha) \le A(n,\alpha)$
we obtain the bounds 
   \begin{equation}\label{eq:angle}
     A(n-1,\omega(\theta,\phi))\le A(n,\theta,[\psi,\phi])\le
        A(n,\theta,\phi)
        \le A(n, \omega(\theta,\phi))
   \end{equation}
which imply this theorem.

Inequality (\ref{eq:angle}) can be rewritten in a somewhat more visual way.
Consider two points $x_1,x_2\in \Cap(n,\phi)$ with angular
distance 
$\theta.$ Let $d$ be the Euclidean distance between them. 
From (\ref{eq:cos}) the angular distance between their images
under $\Pi_n$ is $\theta'\le\omega(\theta,\phi)$ and the Euclidean distance
equals $2\sin\nicefrac{\theta'}2\le d/\sin\phi.$ 
Hence the minimum distance of the image code $\Pi_n(C)\subset S^{n-2}$ 
is at most $d/\sin\phi.$ Denote by $N(n,d, \phi)
=A(n,\arccos(1-\nicefrac{d^2}2),\phi)$ the maximum number
of points in a spherical cap code with minimum Euclidean distance $d$
and let $N(n,d)$ be the same for the sphere. 
An equivalent form of (\ref{eq:angle}) for spherical caps is as follows:
                    $$ N(n-1,d/\sin\phi)\le N(n,d,\phi)\le
                          N(n,d/\sin\phi).$$

A somewhat weaker upper bound on $A(n,\theta,\phi)$
than in (\ref{eq:angle}) is
given by (\ref{eq:ls1}). This bound is nevertheless sufficient
to establish one part of the asymptotic claim of Theorem \ref{cor:eq}.

\section{Codes in hemispheres} \label{sect:hemi}

\begin{figure}[H]\begin{center}\includegraphics[width=2.5in]{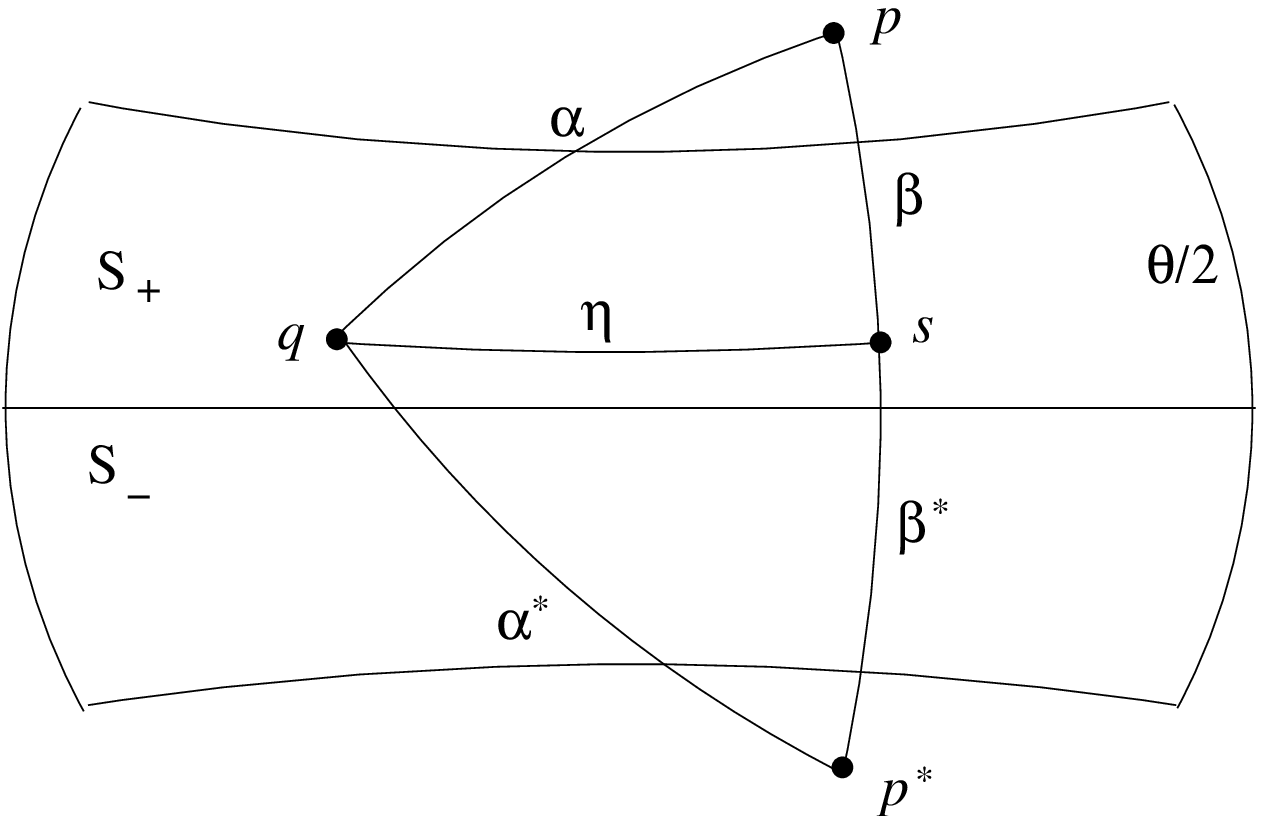}
\end{center}
\caption{}\label{fig:dd}
\end{figure}

Here we consider upper bounds on $B(n,\theta)=A(n,\theta,\pi/2).$
Let $C\subset S_+=\Cap(n,\pi/2)$ be a $\theta$-code in the hemisphere.
Denote 
   \begin{equation}\label{eq:not}
     C([\alpha,\beta])=C\cap Z(n,[\alpha,\beta]), \quad C(\alpha)=
          C\cap \Cap(n,\alpha).
   \end{equation}
\begin{theorem}\label{thm:half} Let $\theta<\nicefrac\pi2, 
\delta=(\pi-\theta)/2$ and let
$C\subset S_+$ be a $\theta$-code. Then
     $$
       |C([\delta,\pi/2])| + 2 |C(\delta)|\le A(n,\theta).
     $$
\end{theorem}
\begin{proof} Let $a=|C([\delta,\pi/2])|, b=|C(\delta)|.$ 
For a point $x=(x_1,\dots,x_{n-1},x_n)$ denote by $x^\ast=
(x_1,\dots,x_{n-1},-x_n)$ its reflection about the equator.
Let $C^\ast(\delta)=\{x^\ast : x\in C(\delta)\}$ be the reflection of the
code $C(\delta)$. Consider the code $Q=C\cup C^\ast(\delta)$.  We claim
that $Q$ is a $\theta$-code. Referring to Fig.~\ref{fig:dd} this amounts to showing that $\dist(q,p^\ast)
\ge\theta$ if $\dist(q,p)\ge\theta.$ To prove this, we choose the point $s$
so that the angle $\angle qsp=90^\circ$ and use (\ref{eq:cos}) as follows:
   $$
    \cos\theta\ge \cos \alpha=\cos\eta\cos\beta\ge\cos\eta\cos\beta^\ast
        =\cos\alpha^\ast.
 $$
This proves that $Q$ is a $\theta$-code. 
Then $a+2b=|Q|\le A(n,\theta).$ 
\end{proof}
\begin{corollary}\label{cor:half} 
   $$
     B(n,\theta)\le \nicefrac12(A(n,\theta,[(\pi-\theta)/2,\pi/2])+
       A(n,\theta))
   $$
\end{corollary}
\begin{proof} Using the notation of the previous theorem, we have
   $$
    2|C|=2a+2b\le a+A(n,\theta)\le A(n,\theta,[\delta,\pi/2])+A(n,\theta).
   $$
\hfill\end{proof}

\begin{corollary}
     \begin{equation}\label{eq:hemi}
        B(n,\theta)\le 
         \frac{A(n-1,\tilde\theta)+A(n,\theta)}{2}, \quad
\cos{\tilde\theta}=\frac{\cos{\theta}}{\cos{\nicefrac{\theta}{2}}}.
      \end{equation}
\end{corollary}
\begin{proof}  
We use Theorem \ref{thm:avz}. 
In our case $\phi=\pi/2, \psi=(\pi-\theta)/2,$ therefore, 
$\theta>\phi-\psi>\theta/2$, so part (a) of this theorem
applies. Substituting the values of $\phi,\psi$ in the inequality
  $  
A(n,\theta,[\psi,\phi])\le A(n-1,\omega(\theta,\phi,\psi)),
  $
we obtain $A(n,\theta,[\psi,\pi/2])\le A(n-1,\tilde\theta).$
The proof is concluded by using this estimate in Corollary \ref{cor:half}.
\end{proof}

Using the value $\theta=60^\circ,$ we obtain estimates on the one-sided
kissing number $B(n).$ In particular, since $A(n,60^\circ,30^\circ)=2,$ 
the last corollary yields
\begin{corollary}\label{cor:ks}
$$B(n)\le \min \left[(1/2)(A(n-1,\eta_0)+k(n)),\; k(n)-2\right], 
\quad \eta_0:=\arccos{\frac{1}{\sqrt{3}}}\approx 54.74^\circ.$$
\end{corollary}
\noindent Clearly, $B(2)=4.$ Let us use this bound for $n=3,4.$
\medskip

$n=3.$ We have $k(3)=12, A(2,\eta_0)=6$, then $B(3)\le
9$. On the other hand, $B(3)\ge 9,$ so
$B(3)=9$. Note that in this case the bound is sharp.

\medskip

 $n=4.$ Recently, K. Bezdek
\cite{bez04,bez05} proved that $B(4)=18$ or 19, and conjectured that
$B(4)=18.$
It was proved, also recently \cite{mus03,mus05b}, that
$k(4)=24$. Delsarte's linear programming method gives $A(3,\eta_0)\le
15$. Thus $B(4)\le 19$. The proof that $B(4)=18$ given 
in \cite{mus05c} is based on an extension of Delsarte's method.

For higher dimensions we can rely on the known bounds for spherical codes.
Denote by $\hat g_n$ an upper bound on $A(n,\eta_0)$ given by 
Delsarte's linear programming method, and by $\hat k_n$ the known upper 
bounds on $k(n)$ (see, e.g., Table 1.5 in \cite{con88}). 
Then Corollary \ref{cor:ks} implies 
    $$
      B(n)\le \frac{\hat g_{n-1}+\hat k_n}2.
    $$
This gives the following bounds:
    $$
     B(5)\le 39,\quad B(6)\le 75,\quad B(7)\le 135,\quad B(8)\le 238.
    $$
For instance, for $n=8$ we have $\hat g_7=236$ (obtained by Delsarte's
method with a polynomial of degree $11$) and $k_8=240.$
Note that $ \hat g_{n-1}>\hat k_n$ for $n>8$; therefore, for these 
dimensions we just have the bound $B(n)\le k(n)-2.$

However, even for $5\le n\le 8$ these bounds are not sharp. Our 
conjectures for $n=5, 8$ are
          $$B(5)=32, \quad B(8)=183.$$

\section{A bound on spherical cap codes}\label{sect:omega}

In this section the methods and results developed above will be used to
derive another bound on spherical cap codes. Given a code 
$C\subset \Cap(n,\phi),$ our plan is to first map the cap on the 
hemisphere $S_+$ and then use the results of the previous section together
with some additional ideas. In the next theorem we use notation (\ref{eq:not}).

\begin{theorem} \label{thm:strf}
Let $0\le \phi\le\pi/2, \omega_1:=\omega(\theta,\phi),
\omega_2:=\omega(\theta,\phi,\psi),$ where 
     $\psi=\phi(1-{\omega_1}/\pi)$.
Then for any $\theta$-code $C\subset \Cap(n,\phi)$ we have
   \begin{equation}\label{eq:p}
     |C([\psi,\phi])|\le A(n-1,\min(\omega_1,\omega_2))
   \end{equation}
and 
   \begin{equation}\label{eq:pp}
    |C([\psi,\phi])|+2|C(\psi)|\le A(n,\omega_1).
   \end{equation}
\end{theorem}
\begin{proof} Consider the mapping $T_s, s=\pi/2\phi$ that sends the cap to 
the hemisphere
$S_+.$ By Theorem \ref{thm:recursion}, the code $C$ is mapped to a 
code $T_s(C) \subset S_+$ with distance $\omega_1.$ 
Now (\ref{eq:pp}) is implied by Theorem \ref{thm:half}. 

To prove (\ref{eq:p}), let us bound above $|C([\psi,\phi])|$ in terms of 
$\theta$ and $\phi.$ Consider the action on the code $C([\psi,\phi])$
of the orthogonal projection $\Pi_n$ on the sphere from $e_n$ on 
the equator. The code $D=\Pi_n(C([\psi,\phi]))$ is a spherical code
in $n-1$ dimensions. Given two points $x,y\in C([\psi,\phi])$ such
that $\dist(x,e_n)=\alpha_1,\dist(y,e_n)=\alpha_2, \dist(x,y)=\beta,$
the distance between their images under $\Pi_n$ is given by
$\omega(\beta,\alpha_1,\alpha_2).$ 
As shown in the proof of Lemma \ref{lemma:lower}, this function is 
monotone on
$\alpha_1$ (or $\alpha_2$) if the other two arguments are fixed, so
its minimum is attained at one of the boundaries. Therefore, the
distance of the code $D$ is determined 
according to one of the following two cases:
\begin{enumerate}\item[(i)] There are two points $c_1,c_2\in C$ lying on
the boundary of the cap and $\theta$ away from each other. Their images in
the code $D$ are two points on the equator at distance $\omega_1.$

\item[(ii)] There are two points $c_1,c_2$ in $C$ at distance $\theta$
such that $T_s(c_1)=x_1, \;T_s(c_2)=x_2$ (see Fig.~\ref{fig:min}). 

\begin{figure}[H]\begin{center}\includegraphics[width=3in]{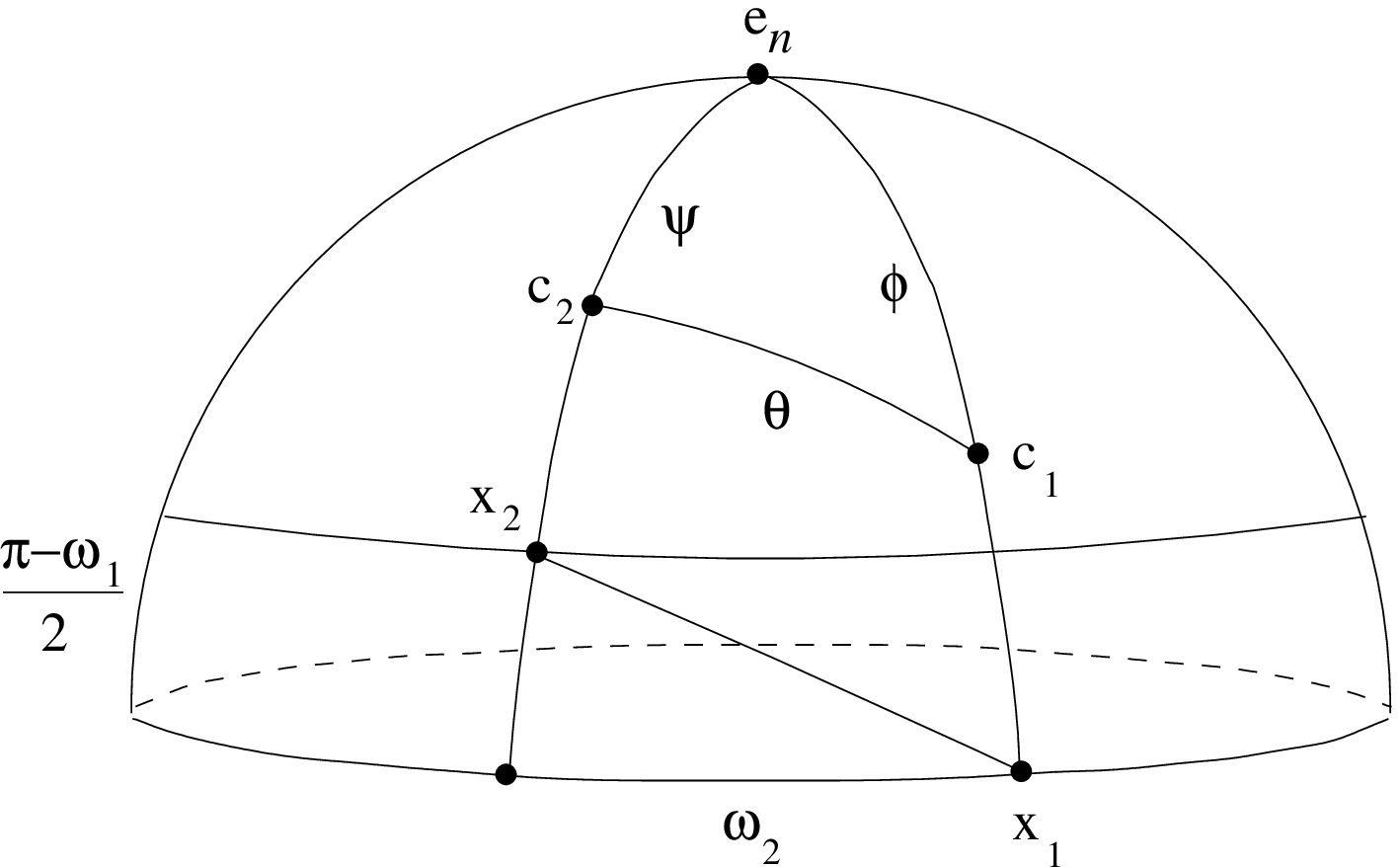}
\end{center}
\caption{}\label{fig:min}
\end{figure}

Upon projecting $x_2$ on the
equator, the distance $\dist(\Pi_n(x_2),x_1)=\omega(\theta,\phi,\psi),$
   where $\psi=\dist(e_n,c_2)$ satisfies
    $$
      \psi s=\frac{\pi-\omega_1}{2},
    $$
i.e., $\dist(\Pi_n(x_2),x_1)=\omega_2.$
\end{enumerate}
Therefore, $D$ is a code with distance at least $\min(\omega_1,\omega_2).$
\end{proof}

\begin{corollary}\label{cor:new}
For all $0<\theta<\phi<\nicefrac \pi 2$
     \begin{equation}\label{eq:new}
      A(n,\theta,\phi) \le \half(A(n-1,\min(\omega_1,\omega_2))+A(n,\omega_1)).
     \end{equation}
\end{corollary}
\begin{proof}
Use (\ref{eq:p}) and (\ref{eq:pp}) to compute
  $$
   2|C([\psi,\phi])|+2|C(\phi)|=|C([\psi,\phi])|+|C(\phi)|+|C(\psi)|\le
           A(n-1,\min(\omega_1,\omega_2))+A(n,\omega_1).
 $$
This gives the claimed result. \end{proof}

\begin{figure}[H]\begin{center}\includegraphics[width=2in]{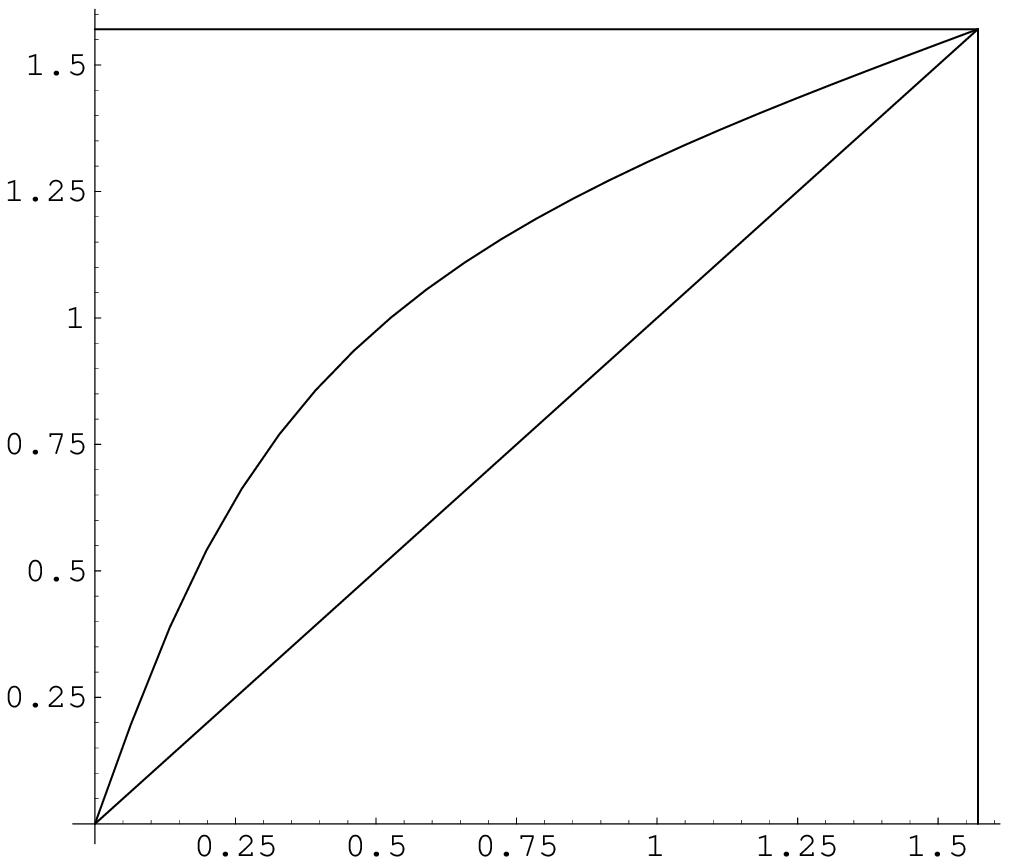}
\put(6,0){{\mbox{\small$\theta$}}}\put(-140,130){{\mbox{\small$\phi$}}}
\put(-50,40){{\mbox{\footnotesize$\theta>\phi$}}}
\put(-87,75){{\mbox{\footnotesize $\omega_1<\omega_2$}}}
\put(-120,100){{\mbox{\footnotesize $\omega_1>\omega_2$}}}
\end{center}
\caption{Area of improvement for bound (\ref{eq:new})}\label{fig:ww}
\end{figure}

Comparison of this bound with the other bounds considered in this paper
is difficult in general. However, it is clear, that it improves upon the
bound (\ref{eq:hemic}) in all cases when 
  $\omega_1\le \omega_2.$ The domain of values of $\theta,\phi$ for which
this holds true is shown in Fig.~\ref{fig:ww}.

{\bf Examples.}
Next we give several examples of the bounds known earlier and obtained in 
this paper, showing that the results of this paper improve the known 
estimate in some range of values of $\theta$ and $\phi$.

\begin{center} 
Bounds on $A(6,\theta,\phi)$ \\[1mm]
\begin{tabular}{|rrrrr|}
\hline &&&&\\
$\theta$\;\;   &$\phi$\;\; &Corollary \ref{cor:agr00} &Theorem \ref{thm:hemi} &Corollary \ref{cor:new} \\[2mm]
\hline &&&&\\[-1mm]     
   $0.02\pi$   &$0.2\pi$ &9069268 &7581886 &4058040\\
     $0.05\pi$  &$0.25\pi$ &145587 &194908 &111241\\
     $0.2\pi$   &$0.4\pi$ &661 &832 &591\\
     $0.25\pi$  &$0.4\pi$ &221 &272 &201\\
     $0.3\pi$   &$0.47\pi$ &174 &138 &115\\[1mm]
\hline
\end{tabular}
\end{center}
\medskip\begin{center}
Bounds on $\log_2 A(n,\theta,\phi)$\\[1mm]
\begin{tabular}{|rrrrrr|}
\hline &&&&&\\
$n$ &$\theta$\;\;   &$\phi$\;\; &Corollary \ref{cor:agr00} &Theorem \ref{thm:hemi} &Corollary \ref{cor:new} \\[2mm]
\hline &&&&&\\[-1mm]     
 20&   $0.05\pi$   &$0.4\pi$ &66.1912 &65.7786 &65.3665\\
 20&    $0.08\pi$  &$0.4\pi$ &52.7212 &52.891 & 52.6557\\
 30&    $0.05\pi$   &$0.4\pi$ &97.75 &97.7008 &97.8523\\
40&    $0.05\pi$  &$0.45\pi$ &132.516 &131.447 &133.45\\
100&   $0.05\pi$  &$0.45\pi$ &322.205 &321.917 &333.28\\[1mm]
\hline
\end{tabular}
\end{center}

\remove{ Let us show that the bound of Theorem \ref{thm:hemi}
is sometimes better than the result of Corollary \ref{cor:agr00} \cite{agr00}.

We wish to estimate $\log_2 A(n,\theta,\phi).$
Let $n=40, \theta=13^\circ, \phi=80^\circ.$ From (\ref{eq:capsavz}) 
we obtain the estimate $111.96.$ For the same parameters, 
(\ref{eq:hemic}) gives $110.59.$ Letting 
$n=257, \theta=20^\circ, \phi=85^\circ,$ we obtain $523.554$ 
and $522.95$ from (\ref{eq:capsavz}) and (\ref{eq:hemic}), respectively.
Observe also that the bound (\ref{eq:hemic}) is notably easier to compute
than (\ref{eq:capsavz}) which is generally a multi-step recursive 
calculation.}

\medskip 
To compute the bounds in examples one needs to use some upper bound on 
$A(n,\theta).$
We have used the bound of Levenshtein \cite{lev83a},\,\cite[p.618]{lev98}
which is the best known universal upper bound on spherical codes.
\section{A general method of bounding the size of cap codes}\label{sect:lp}
In this section we generalize the method of Sections \ref{sect:hemi}, 
\ref{sect:omega} to develop a general approach
to bounding the size of codes in spherical caps. 

Let $C\in \Cap(\theta,\phi)$ be a $\theta$-code. 
Let $0=\phi_0<\phi_1<\dots<\phi_k<\phi_{k+1}=\phi.$

To bound above the size of the code $C$ we partition the cap as follows:
  $$
    \Cap(n,\phi)=\Cap(n,\phi_1)\cup Z(n,[\phi_1,\phi_2])\cup
     \dots\cup Z(n,[\phi_k,\phi_{k+1}]).
  $$
Let $p_i:=|C \cap Z(n,[\phi_i,\phi_{i+1}])|$ be the size of the code in the
strip which can be estimated by the methods of \cite{agr00} and this paper.
In particular, denote by $a_{i,j}$ an upper bound on the size of the
code in the strip $Z(n,[\phi_i,\phi_j])$ and let $a_i:=a_{i,i+1}.$
Suppose that we can compute these estimates for some subset $P$ of pairs
$(i,j).$ 
We have the following linear constraints.
  \begin{equation}\label{eq:c1}
     0\le p_i\le a_i,\quad (i,i+1)\in P; \qquad p_i+\dots p_{j-1}\le a_{i,j},
  \quad(i,j)\in P.
  \end{equation}
Another set of linear inequalities can be obtained from the arguments
similar to the proof of Theorems \ref{thm:half}, \ref{thm:strf}. Namely, let 
$r(i)=\max\{j: s\phi_j\le (\pi-\omega(\theta,\phi_{i+1}))/2\}.$ Consider the
codes $C'=C(\phi_{r(i)+1})$ and $C''=C([\phi_{r(i)+1},\phi_{i+1}]).$
Upon stretching the cap to the hemisphere by applying the mapping $T_s$
and using (\ref{eq:pp}), we obtain $2|C'|+|C''|\le b_i,$ where
$b_i$ is the size of spherical code with distance 
$\omega(\theta,\phi_{r(i)+1}).$ We then obtain inequalities
  \begin{equation}\label{eq:c2}
    2(p_0+p_1+\dots+p_{r(i)})+p_{r(i)+1}+\dots+p_i\le b_i, \quad i\in\Phi,
  \end{equation}
where $\Phi$ is a subset of indices for which we perform the described
procedure.

We summarize the arguments of this section in the following
\begin{theorem} The size of the code $C$ is bounded above by the
solution of the linear programming problem
  $$
     p_0+p_1+\dots+p_k\to\max
  $$
under the constraints (\ref{eq:c1}), (\ref{eq:c2}).
\end{theorem}

Note that Theorem \ref{thm:strf} constitutes a solution of this problem
in the case $k=0.$

The set $\{\phi_i\}$ can be optimized in computations. It can be formed,
for instance, by taking all the angles $\gamma$ produced by the recurrent
calculation in Corollary \ref{cor:agr00} together with some additional
breakpoint angles. 

\section{Applications}
We have discussed the applications of the bounds on codes in caps
to the kissing number and the one-sided kissing number problems.
As remarked in the introduction,
the methods developed in this paper can be also
useful in the problems of estimating
the packing density in $\reals^n$ and of the size of constant weight codes.
We end the paper with brief remarks on these applications.

\subsection{Spherical codes and packing density}
Let $\Delta_n$ be the density of packing the $n$-dimensional real space
with equal nonoverlapping balls. A classical problem in geometry is
to compute $\Delta_n$ for a given $n$ and for $n\to\infty.$
It is known \cite[p.265]{con88} that
      \begin{equation}\label{eq:delta}
     \Delta_{n}\le \nicefrac 12(\sin\nicefrac\theta2)^{n}A(n+1,\theta) \quad
           (0<\theta\le\pi).
      \end{equation}
Using Lemma \ref{lemma:be} we now obtain
the estimate
    \begin{equation}\label{eq:deltaBE}
     \Delta_{n}\le \nicefrac 12\;
          (\sin\nicefrac\theta2)^{n}\Omega_{n+1} \frac{A(n+1,\theta,\phi)}{\Omega_{n+1}
         (\phi)}\quad
           (0<\theta\le\pi, \theta/2\le\phi\le \pi/2).
    \end{equation}
To compute upper bounds on $\Delta_n$ for a given value of $n$ we
have a choice of using (\ref{eq:delta}) with the known bounds
on $A(n+1,\theta),$ or using (\ref{eq:deltaBE}) together with the bounds
on spherical cap codes considered in this paper. Observe that the
best known estimate of the packing density for $n\to\infty$ \cite{kab78} 
is obtained
by employing (\ref{eq:deltaBE}) together with the
bound (\ref{eq:ls1}) on cap codes. Therefore, inequality (\ref{eq:deltaBE})
coupled with better bounds on cap codes derived in this paper
will also improve the density estimates for finite (but possibly large)
values of $n$.

\subsection{Constant weight codes} A constant weight binary code is a subset
of $\{0,1\}^n$ formed of vectors with a fixed number, say $w$, of ones.
Denote by $A(n,d,w)$ the maximum size of a constant weight code with
minimum Hamming distance $d.$ Computing or estimating the numbers
$A(n,d,w)$ is a problem with a long history in coding theory, summarized
in \cite{agr00}. The most studied region of parameters is $n\le 28$ for
which the most recent tables are published in \cite{agr00}.
For larger $n$ bounds on $A(n,d,w)$ are tabulated in \cite{smi06}. 
In \cite{agr00} the problem of bounding
above $A(n,d,w)$ was reduced to bounds on spherical cap codes which led to
several improvements of the tables for short lengths. We intend to use
the new bounds on cap codes derived in this paper to further improve the
tables.

As a final remark, note that all the bounds on cap codes considered in 
this paper
rely on bounds on the maximum size $A(n,\theta)$ of spherical codes.
In calculations, apart from the Levenshtein bound and related results
it is possible to use a direct solution of Delsarte's linear programming
problem relying on a method developed in \cite{odl79}.

\providecommand{\bysame}{\leavevmode\hbox to3em{\hrulefill}\thinspace}


\begin{thebibliography}{10}

\bibitem{agr00}
E.~Agrell, A.~Vardy, and K.~Zeger, \emph{Upper bounds for constant-weight
  codes}, IEEE Trans. Inform. Theory \textbf{46} (2000), 2373--2395.

\bibitem{bez04}
K.~Bezdek, \emph{Sphere packings in 3-space}, Proceedings of the COE Workshop
  on Sphere Packings, Kyushu University Press, 2004, pp.~32--39.

\bibitem{bez05}
\bysame, \emph{Sphere packing revisited}, European J. Comb. \textbf{27} (2006),
  no.~6, 864--883.

\bibitem{boy94}
P.~Boyvalenkov, \emph{Small improvements of the upper bounds of the kissing
  numbers in dimensions {$19$}, {$21$} and {$23$}}, Atti Sem. Mat. Fis. Univ.
  Modena \textbf{42} (1994), no.~1, 159--163.

\bibitem{con88}
J.~H. Conway and N.~J.~A. Sloane, \emph{Sphere packings, lattices and groups},
  Springer-Verlag, New York-Berlin, 1988.

\bibitem{dan86}
L.~Danzer, \emph{Finite point-sets on {${\bf S}^2$} with minimum distance as
  large as possible}, Discrete Math. \textbf{60} (1986), 3--66.

\bibitem{del73}
P.~Delsarte, \emph{An algebraic approach to the association schemes of coding
  theory}, Philips Research Repts Suppl. \textbf{10} (1973), 1--97.

\bibitem{del77b}
P.~Delsarte, J.~M. Goethals, and J.~J. Seidel, \emph{Spherical codes and
  designs}, Geometriae Dedicata \textbf{6} (1977), 363--388.

\bibitem{fej53}
L.~Fejes~{T}{\' o}th, \emph{Lagerungen in der {E}bene, auf der {K}ugel und in
  {R}aum}, Springer-Verlag, 1953.

\bibitem{kab78}
G.~Kabatyansky and V.~I. Levenshtein, \emph{Bounds for packings on the sphere
  and in the space}, Problems of Information Transmission \textbf{14} (1978),
  no.~1, 3--25.

\bibitem{lev75}
V.~I. Levenshtein, \emph{The maximal density of filling an $n$-dimensinal
  {E}uclidean space with equal spheres}, Mat. Zametki \textbf{18} (1975),
  no.~2, 301--311, English translation in Math. Notes 18 (1975), no. 1--2,
  765--771.

\bibitem{lev83a}
\bysame, \emph{Bounds for packings of metric spaces and some of their
  applications}, Problemy Kibernet. (1983), no.~40, 43--110 (In Russian).

\bibitem{lev98}
\bysame, \emph{Universal bounds for codes and designs}, Handbook of Coding
  Theory (V.~Pless and W.~C. Huffman, eds.), vol.~1, Elsevier Science,
  Amsterdam, 1998, pp.~499--648.

\bibitem{mus05b}
O.~R. Musin, \emph{The kissing number in four dimensions}, preprint, April
  2005, arxiv.org/math.MG/0309430.

\bibitem{mus05c}
\bysame, \emph{The one-sided kissing number in four dimensions}, preprint,
  November 2005, arxiv.org/math.MG/0511071.

\bibitem{mus03}
\bysame, \emph{The problem of the twenty-five spheres}, Russian Math. Surveys
  \textbf{58} (2003), 794--795.

\bibitem{mus04}
\bysame, \emph{An extension of {D}elsarte's method. {T}he kissing number
  problem in three and four dimensions}, Proceedings of the COE Workshop on
  Sphere Packings, Kyushu University Press, 2004, pp.~1--25.

\bibitem{mus05a}
\bysame, \emph{The kissing problem in three dimensions}, Discrete and
  Computational Geometry \textbf{35} (2006), 375--384.

\bibitem{odl79}
A.~M. Odlyzko and N.~J.~A. Sloane, \emph{New bounds on the number of spheres
  that can touch a unit sphere in $n$ dimensions},  \textbf{26} (1979),
  210--214.

\bibitem{pfe05}
F.~Pfender, \emph{Improved {D}elsarte bounds via extension of the function
  space}, preprint, 2005, arxiv.org/math.CO/0501493.

\bibitem{ran55}
R.~A. Rankin, \emph{The closest packing of spherical caps in $n$ dimensions},
  Proc. Glasgow Math. Assoc. \textbf{2} (1955), 139--144.

\bibitem{sch51}
K.~Sch{\"u}tte and B.~L. van~der Waerden, \emph{Auf welcher {K}ugel haben
  {$5$}, {$6$}, {$7$}, {$8$} oder {$9$} {P}unkte mit {M}indestabstand {E}ins
  {P}latz?}, Math. Ann. \textbf{123} (1951), 96--124.

\bibitem{sch53}
\bysame, \emph{Das {P}roblem der dreizehn {K}ugeln}, Math. Ann. \textbf{125}
  (1953), 325--334.

\bibitem{sid74}
V.~M. Sidelnikov, \emph{New estimates for the closest packing of spheres in
  $n$-dimensional {E}uclidean space}, Mat. Sb. (N.S.) \textbf{95(137)} (1974),
  148--158, 160.

\bibitem{smi06}
D.~H. Smith, L.~A. Hughes, and S.~Perkins, \emph{A new table of constant weight
  codes of length greater than 28}, Electron. J. Combin. \textbf{13} (2006),
  Article A2.

\end{thebibliography}
\end{document}